\documentclass[intlim,righttag,11pt]{amsart}
\usepackage{stmaryrd}
\usepackage{amscd}
\usepackage{amssymb}
\usepackage[all]{xy}
\RequirePackage{mathrsfs}
\def\jcdot{{\scriptscriptstyle\bullet}}

\def\invlim{\mathop{\vtop{\ialign{##\crcr$\hfill{\lim}\hfil$\crcr
\noalign{\kern1pt\nointerlineskip}\leftarrowfill\crcr\noalign
{\kern -3pt}}}}\limits}
\def\dirlim{\mathop{\vtop{\ialign{##\crcr$\hfill{\lim}\hfil$\crcr
\noalign{\kern1pt\nointerlineskip}\rightarrowfill\crcr\noalign
{\kern -3pt}}}}\limits}

 \def\verylomapr#1{\smash{\mathop{\relbar\joinrel\relbar\joinrel\relbar\joinrel\longrightarrow}\limits^{#1}}}

\def\phi{\varphi}
\def\epsilon{\varepsilon}
\let\mathcal\mathscr
\usepackage{fge}
\usepackage{bbm}
\renewcommand{\setminus}{\mathbin{\fgebackslash}}

\newtheorem{theorem}[equation]{Theorem}

 \newtheorem{corollary}[equation]{Corollary}
\newtheorem{conjecture}[equation]{Conjecture}

\theoremstyle{definition}

\theoremstyle{remark}
\newtheorem{remark}[equation]{Remark}

\newtheorem{example}[equation]{Example}
\newtheorem*{acknowledgments}{Acknowledgments}

\newcommand{\what}{\widehat}

\renewcommand{\phi}{\varphi}

\newcommand{\pst}{\rm pst}

\newcommand{\ovk}{\overline{K} }

 \newcommand{\cont}{\operatorname{cont}  } 
 
  \newcommand{\proeet}{\operatorname{pro\acute{e}t}  }

 \newcommand{\eet}{\operatorname{\acute{e}t} }

 \newcommand{\dlog}{\operatorname{dlog} }

 \newcommand{\nr}{\operatorname{nr} }

 \newcommand{\Ext}{\operatorname{Ext} }
 \newcommand{\Rep}{\operatorname{Rep} }

 \newcommand{\Gal}{\operatorname{Gal} }
 
 \newcommand{\can}{ \operatorname{can} }

 \newcommand{\id}{ \operatorname{Id} }

\newcommand{\st}{\operatorname{st} }

\newcommand{\hk}{\operatorname{HK} } 
\newcommand{\dr}{\operatorname{dR} }

 \newcommand{\crr}{\operatorname{cr} }
 \newcommand{\gr}{\operatorname{gr} }

  \newcommand{\htt}{\operatorname{HT}}

 \newcommand{\scc}{{\mathcal{C}}}

 \newcommand{\so}{{\mathcal O}}

 \newcommand{\wt}{\widetilde}
 \newcommand{\wh}{\widehat}

  \newcommand{\C}{ {\mathbf C} }
 \newcommand{\Z}{ {\mathbf Z} }
    
   \newcommand{\Q}{ {\mathbf Q}}
   \newcommand{\N}{{\mathbf N}}
         
   \def\B{{\bf B}}
      \def\A{{\bf A}}

\numberwithin{equation}{section}

\begin{document}
 %\title[Towards $p$-adic Hodge Theory of analytic varieties: a survey]{Towards $p$-adic Hodge Theory for analytic varieties: a survey}
 \title[Hodge Theory of $p$-adic varieties]{Hodge Theory of $p$-adic varieties: a survey}
  \author{Wies{\l}awa Nizio{\l}}
\address{CNRS, IMJ-PRG, Sorbonne Universit\'e, 4 place Jussieu, 75005 Paris, France}
\email{wieslawa.niziol@imj-prg.fr}
\thanks{The authors' research was supported in part by the grant ANR-19-CE40-0015-02 COLOSS} 
 \date{\today}

\begin{abstract}$p$-adic Hodge Theory is one of the most powerful tools  in modern Arithmetic Geometry. In this survey, 
we will review $p$-adic Hodge Theory for algebraic varieties,  present  current developments in $p$-adic Hodge Theory for analytic varieties, and discuss some of its applications to problems in Number Theory. 
This is an extended version of a talk at the Jubilee Congress for the 100th anniversary of the Polish Mathematical Society, Krak\'ow, 2019.
\end{abstract}
 \setcounter{tocdepth}{2}

 \maketitle
 \tableofcontents
 \section{Introduction}$p$-adic Hodge Theory is one of the most powerful tools  in modern Arithmetic Geometry. In this survey, 
 we will review $p$-adic Hodge Theory for algebraic varieties,  present  current developments in $p$-adic Hodge Theory for analytic varieties, and discuss some of its applications to problems in Number Theory.

 This is an extended version of a talk at the Jubilee Congress for the 100th anniversary of the Polish Mathematical Society, Krak\'ow, 2019.

%$p$-adic Hodge Theory is one of the most powerful tools we have currently in Arithmetic Geometry. We will discuss some of its applications.
\begin{acknowledgments} We would like to thank the organizers of the Congress for giving us the opportunity to present this talk. Special thanks go to 
Pierre Colmez for many wonderful discussions concerning the content of this paper. 
\end{acknowledgments}
 \section{Algebraic varieties}
 We will first consider algebraic varieties, starting with the classical situation over complex numbers. 
 \subsection{A comparison theorem over $\C$} 
Let $X$ be a smooth and projective algebraic variety over the rational numbers $\Q$.
 Recall that, by the classical de Rham theorem,
 there exists a nondegenerate pairing
  $$
  H^n_{\dr}(X_\C)\times H_n(X(\C),\C)\to \C,\quad (\omega,\gamma)\mapsto \int_{\gamma}\omega.
  $$
  Here $H_n(X(\C),\C)$ is the  singular homology with complex coefficients $\C$ of the topological space $X(\C)$ and    $H^n_{\dr}(X_{\C}) $ is the algebraic de Rham cohomology of $X_{\C}$:
  $$
  H^n_{\dr}(X_{\C}) :=H^n(X_{\C},\so_{X_{\C}}\to \Omega^1_{X_{\C}}\to\Omega^2_{X_{\C}}\to\cdots).
$$
The pairing is obtained by integrating differential forms along cycles. Using resolutions of singularities, and changing de Rham cohomology and singular homology appropriately,  this pairing can be extended to all algebraic varieties. So obtained values are (classically) called  {\em periods}. 
Hence, we can say that:
\vskip.2cm
\begin{center}\fbox{
 $\C$ {\it  contains all periods of algebraic varieties over $\Q$}. 
 }
 \end{center}
\vskip.2cm

\begin{example} Here are a couple of  examples of  periods: 
 \begin{enumerate}
 \item Let $\gamma$ denote the unit circle in the complex plane. We have
$$\int_{\gamma}\frac{{\rm d}z}{z}=2\pi i.$$
To put this integral in the  context discussed above, take $X={\mathbb G}_{m,\Q}={\mathbb A}^1_{\Q}\setminus \{0\}$, the punctured affine space. Then $\gamma\in H_1(X(\C),\C)$, $\omega={\rm d}z/z\in H^1_{\dr}(X_{\C})$, and $(\omega,\gamma)=2\pi i$.
\item Consider now the integral in the complex plane 
$$
2\int_1^{+\infty}\frac{dx}{\sqrt{x^3-x}} =\frac{\Gamma(1/4)\Gamma(1/2)}{\Gamma(3/4)},$$
where $\Gamma(-)$ denotes Euler's $\Gamma$-function.
This is a period of the elliptic curve $E$ with equation $y^2=x^3-x$. (Modulo a sign) this is the integral of $\omega={\rm d}x/y$ along the  path in $E(\C)$ whose projection in the projective space ${\mathbb P}^1(\C)=\C\cup\{\infty\}$ (by  the map $(x,y)\mapsto x$)    consists of the segment $(\infty, 1+\varepsilon]$, followed by the  circle of center $1$ and radius $\varepsilon$, followed by the segment $[1+\varepsilon, \infty)$. This integral does not depend on $\varepsilon$  and when $\varepsilon \to 0$ the contribution from the circle tends towards $0$. Since $\sqrt{x^3-x}$ changes the sign going around the circle with center  $1$,  the integrals $\int^{\infty}_{1}$ and $\int^{1}_{\infty}$ are equal. Hence   the integral along the whole path in $E(\C)$  is $2\int_1^{+\infty}\frac{dx}{\sqrt{x^3-x}} $.
\end{enumerate}
%\begin{align*}
%& \mbox{Dually :} \quad \quad H^n_{\dr}(X)\otimes_K\C\simeq  H^n_{B}(X(\C),\Q)\otimes_{\Q}\C
 %\end{align*}

 We note that we started here with an algebraic variety $X$ over $\Q$, then we completed $\Q$ with respect to its natural Archimedean valuation, i.e.~the usual absolute value,
 to obtain the field of real numbers $\mathbf{R}$ and the base changed variety 
 $X_{\mathbf{R}}$, and then we passed to the algebraic closure $\C$ of $\mathbf{R}$ and the corresponding variety $X_{\C}$. We can illustrate this process by the picture:
 $$
 \Q\mapsto \wh{\Q}\simeq \mathbf{R}\hookrightarrow\overline{\mathbf{R}}\simeq \C.
 $$
 \end{example}
  \subsection{\'Etale cohomology and associated Galois representations}
Grothendieck defined \'etale cohomology as an algebraic replacement of singular cohomology. 
Before we proceed to review its properties let us make a small digression. 
  \subsubsection{Digression: non-Archimedean completions}
Besides the Archimedean valuation on $\Q$,  we also have non-Archimedean valuations\footnote{By Ostrowski's theorem these valuations, taken together with the natural Archimedean valuation, constitute  all nontrivial valuations on $\Q$ (taken up to equivalence).}   indexed by  prime numbers. The relevant picture is now:
 $$
 \Q\mapsto \what{\Q}=: \Q_p\hookrightarrow\overline{\Q}_p\hookrightarrow \wh{\overline{\Q}_p}=:\C_p.
 $$
Here: 
 \begin{enumerate}
 \item $p$ is a prime number. If $|\jcdot|=|\jcdot|_p$ is the  $p$-adic norm on $\Q$, normalized with $|p|=p^{-1}$, we  have $|xy|=|x||y|$ and $|x+y|\leq\max (|x|,|y|)$.
\item  $\Q_p$ is the completion of $\Q$ for the $p$-adic norm $|\jcdot|$ (introduced by Hensel, 1897). We have:
  \begin{align*}
 & \Z_p:=\{x\in\Q_p| |x|\leq 1\},\quad \Z_p\simeq \invlim_n\Z/p^n,\\
& \Z_p"="\{0,1,\ldots,p-1\}[[p]],\\
 & \Q_p=\Z_p[1/p],\quad x\in\Q_p\Rightarrow  x=\sum_{n\geq n_0}x_np^n, \,x_n\in\{0,\ldots,p-1\}.
  \end{align*}
  \item
 $\overline{\Q}_p$ is the algebraic closure of $\Q_p$; the norm $|\jcdot|$ extends uniquely to $\overline{\Q}_p$; the Galois group $G_{\Q_p}:=\Gal(\overline{\Q}_p/\Q_p)$ acts via isometries. 
  $\overline{\Q}_p$ is infinite dimensional (the polynomial $x^n-p$ is irreducible in $\Q_p[x]$),
and  not complete for $|\jcdot|$.
  \item  $\C_p$ is the completion of $\overline{\Q}_p$ and we have  
$G_{\Q_p}={\rm Aut}_{\cont}(\C_p)$. $\dim_{\Q_p}\C_p$ is not countable. 
Axiom of choice produces isomorphisms of abstract fields $\C_p\simeq \C$.
\end{enumerate}
\subsubsection{\'Etale cohomology}
  Let us now come back to  the nondegenerate pairing: 
  $$
  H^n_{\dr}(X_\C)\times H_n(X(\C),\C)\to \C,\quad (\omega,\gamma)\mapsto \int_{\gamma}\omega.
  $$
  Using singular (Betti) cohomology, which is dual to singular homology, this can be written as a $\C$-linear  isomorphism: 
  \begin{equation}
  \label{drpair}
     H^n_{\dr}(X)\otimes_{\Q_p}\C\simeq  H^n_{B}(X(\C),\Q)\otimes_{\Q}\C.
  \end{equation}
  
     Fix a prime $p$. Recall that, for $p$-adic coefficients $\Q_p$, we have the isomorphism 
  \begin{equation}
  \label{betti-et}
  H^n_{B}(X(\C),\Q)\otimes_{\Q}\Q_p\simeq H^n_{\eet}(X_{\overline{\Q}_p},\Q_p)
  \end{equation}
  of Betti  cohomology with  (geometric) \'etale cohomology. 
 The latter cohomology has the following properties: 
  \begin{enumerate}
  \item  The vector space $H^n_{\eet}(X_{\overline{\Q}_p},\Q_p)$ is of finite rank over $\Q_p$;     it inherits a continuous action of $G_{\Q_p}$ (from its natural action on 
$X_{\overline{\Q}_p}$). The $\Q_p$-rank of 
  $H^n_{\eet}(X_{\overline{\Q}_p},\Q_p)$ is equal to the $\Q$-rank of Betti cohomology, hence also, to the $\Q$-rank of de Rham cohomology (see (\ref{betti-et}) and (\ref{drpair})).

    \item Locally in the Zariski topology, \'etale cohomology is the group cohomology of the fundamental group: $H^n_{\eet}(X_{\overline{\Q}_p},\Q_p) \simeq H^n(\pi(X_{\overline{\Q}_p}),\Q_p)$, where $\pi(X_{\overline{\Q}_p})$ is the  algebraic fundamental group (it is the 
profinite completion of the classical fundamental group).
     \end{enumerate}
 \subsubsection{Galois representations coming from \'etale cohomology}
 The Galois action on \'etale  cohomology $H^n_{\eet}(X_{\overline{\Q}_p},\Q_p)$ carries information about:
  \begin{enumerate}
  \item finite extensions of $\Q_p$,
  \item the arithmetic of $X$, for example, its rational points $X(\Q)$ (see Section \ref{applications}).
  \end{enumerate}

 Let us look at some examples of such representations.
  \begin{example}{\em Galois representations on $H^n_{\eet}(X_{\overline{\Q}_p},\Q_p)$}:
  
  (1) {\em Tate twists}: we have the 
cyclotomic character 
$$\chi: G_{\Q_p}\to \Z_p^*: \quad \sigma(e^{\frac{2\pi i}{p^n}})=e^{\chi(\sigma)\frac{2\pi i}{p^n}}, n\geq 1.$$
If $i\in\Z$, denote by $\Q_p(r)$ the $i$'th Tate twist: it is $\Q_p$ with the action of $G_{\Q_p}$ given by  $\chi^r$.
We have
 $$\Q_p(1)\simeq \Q_p\otimes_{\Z_p}
\invlim_n\mathbb{G}_m(\overline{\Q}_p)_{p^n} ,\quad
   \Q_p(1)\simeq H^2_{\eet}(\mathbb{P}^1_{\overline{\Q}_p},\Q_p)^*.$$
Here, the subscript $p^n$ refers to the $p^n$-torsion elements and ${\mathbb P}^1$ denotes the projective space.
  
  (2)  {\em Tate modules}: Let $E$ be an elliptic curve over $\Q$. We have its integral and rational Tate module
  $$
  T_pE:=\invlim_n E(\overline{\Q}_p)_{p^n},\quad V_pE:=\Q_p\otimes_{\Z_p}T_pE, ,\quad \dim _{\Q_p}V_pE=2.
  $$
Cohomologically $$V_pE\simeq H^1_{\eet}(E_{\overline{\Q}_p},\Q_p)^*.$$
  \end{example}
 \subsection{The main question of geometric $p$-adic Hodge Theory} Let us pass  now to the local setting. That is,  $X$ will be now a variety over a finite field extension\footnote{What we will describe holds more generally for complete discrete valuation fields of mixed characteristic with  perfect residue fields.} $K$ of $\Q_p$. The main question of the classical geometric $p$-adic Hodge Theory was the following: \\ 
 
  Does there exist  a $p$-adic period ring $\B$ that contains all periods of algebraic varieties $X$ over $K$, for all $[K:\Q_p]<\infty,$ and a pairing $(\omega,\gamma)\mapsto (\int_{\gamma}\omega)\in \B$ such that
  \begin{enumerate}
  \item $      H^n_{\dr}(X)\otimes_{K}\B\simeq  H^n_{\eet}(X_{\overline{\Q}_p},\Q_p) \otimes_{\Q_p}\B$, 
\item using  the period isomorphism (1) we can recover the Galois representation on the \'etale cohomology $H^n_{\eet}(X_{\overline{\Q}_p},\Q_p)$ from the de Rham cohomology $H^n_{\dr}(X)$.
\end{enumerate}
\begin{remark}It was observed by Tate (circa 1966) that  $\B$ can not be $\C_p$ because the latter   does not contain a $p$-adic analog $(2\pi i)_p$ of $2\pi i$. 
\begin{enumerate}
\item  Naively, we would like to have $(2\pi i)_p=p^n\log e^{\frac{2\pi i}{p^n}}$, $n\geq 1$. But 
the $p$-adic logarithm defined on the the open unit ball $B(1,1^{-})$ with center $1$ and radius $1$,
by the usual formula
$ \log (x)=\sum_{n\geq 1}\frac{(-1)^{n-1}}{n}(x-1)^n $,
is a group homomorphism (as opposed  to its complex analog). It follows that 
we have $\log e^{\frac{2\pi i}{p^n}}=0$ and
an exact sequence
\begin{equation}
\label{logparis}
0\to \mu_{p^{\infty}}\to B(1,1^{-})\stackrel{\log}{\to} \C_p\to 0.
\end{equation}
\item Perhaps less naively, 
$(2\pi i)_p$ should be a period of ${\mathbb G}_{m,\overline{\Q}_p}$ and as such tranform under the Galois action by 
$
\sigma((2\pi i)_p)=\chi(\sigma)(2\pi i)_p,\quad\forall \sigma\in G_{K}.
$
But we have 
$$
\{x\in \C_p\mid \sigma(x)=\chi(\sigma)x,\, \forall \sigma\in G_{K}\}=0,
$$
\end{enumerate}
which is a special case of the following fundamental theorem of Tate: 
\begin{theorem}{\rm(Tate, \cite{Ta})} Let $k\in\Z$. Then 
$$\{x\in \C_p\mid \sigma(x)=\chi(\sigma)^kx, \forall \sigma\in G_{K}\}=\begin{cases} \Q_p & \mbox{ if }k=0,\\
 0 & \mbox{ if }k\neq 0.
\end{cases}
$$
\end{theorem}
\end{remark}
\vskip.2cm
\begin{center}\fbox{
 $\C_p$ {\it  does not contain $(2\pi i)_p$}. 
 }
 \end{center}
\vskip.2cm

 \subsubsection{History of  $p$-adic Hodge Theory for algebraic varieties: the beginnings.}
 \begin{enumerate}
\item  1958-65: Grothendieck defines \'etale cohomology (as an analog of singular cohomology), algebraic de Rham cohomology, and, together with Berthelot,  its refinement: crystalline cohomology.
\item 1967-70: Tate and Grothendieck discover  that, for  an elliptic curve over $K$ (or, more generally an abelian variety), 
its de Rham cohomology and its $p$-adic \'etale cohomology both determine and are determined by its $p^\infty$-torsion (more precisely its $p$-divisible group).
\item 1970: Grothendieck asks
whether there exists an abstract  "mysterious period functor" 
relating directly
$p$-adic \'etale cohomology and de Rham cohomology.
 \item 1979-87: Fontaine constructs complicated period rings and formulates precise conjectures, which are now theorems, concerning this "mysterious functor''. 
 \end{enumerate}
  \subsection{The first approximation: the de Rham period ring $\B_{\dr}$}
  We start with the simplest of $p$-adic  period rings. 
  \subsubsection{The de Rham period ring}
Fontaine (circa 1980) constructed a ring
$
\B^+_{\dr}$, the de Rham period ring, that he conjectures should contain all periods of $p$-adic algebraic varieties. 
In $\B^+_{\dr}$ he distinguished an element $t$ (for Tate) that will become  the $p$-adic analog of $2\pi i$. In particular, the Galois group acts on $\B^+_{\dr}$ and does act on $t$ via the cyclotomic character: $\sigma(t)=\chi(\sigma)t, \sigma\in G_{K}.$

 We list the following properties of $\B^+_{\dr}$:
\begin{enumerate}
\item $\B^+_{\dr}\simeq \C_p[[t]]$ but not in any reasonable way\footnote{There exist  $K$-linear continuous sections $\C_p\to \B^+_{\dr}$ of the surjective map $\theta:  \B^+_{\dr}\to \C_p$ but none  preserving the ring structure. Axiom of choice gives existence of sections preserving the ring structure but they can not be continuous.};  we do, however,  have a short exact sequence
$$
0\to t\B^+_{\dr}\to \B^+_{\dr}\stackrel{\theta}{\to}\C_p\to 0.
$$
\item $\B^+_{\dr}$ is equipped with a descending filtration by the powers of $t$: 
$$\B^+_{\dr}\supset F^n\B^+_{\dr}:=(t^n),\quad  \gr_F^n\B^+_{\dr}\simeq \C_p(n).$$

\item $\B^+_{\dr}$ is  a completion of $\overline{\Q}_p$ involving ``higher derivatives'' \cite{CBdr}.
\end{enumerate}
\begin{example} \label{derham}We will illustrate (3).
Define the norm 
\begin{center}
$
x\in \overline{\Q}_p, \quad |x|_{p,1}:=\sup (|x|_p,|\frac{dx}{dp}|_p).
$
\end{center}
It is a submultiplicative norm : $|xy|_{p,1}\leq |x|_{p,1}|y|_{p,1}$. 
To define $|\frac{dx}{dp}|_p$,
choose a uniformizer $\pi$ of $\Q_p(x)$, and
write $x=Q(\pi)$, for a polynomial $Q\in\Q_p(\mu_N)[X]$, $(N,p)=1$, of minimal degree.  
Let $P\in \Q_p(\mu_N)[X]$ be the
minimal polynomial of $\pi$ (it is an Eisenstein polynomial). Set
$ \frac{dx}{dp}:=\frac{-Q^{\prime}(\pi)}{P^{\prime}(\pi)} $.
Then $\frac{dx}{dp}$ depends on the choices made but not too much, and
$|\frac{dx}{dp}|_p$ does not depend on the choices made.

   The ring $\B^+_{\dr}/t^2$ is the completion of $\overline{\Q}_p$ for the norm $|\jcdot|_{p,1}$. In particular,  $\overline{\Q}_p$ is dense in $\B^+_{\dr}/t^2$. Hence $\B^+_{\dr}/t^2$ is not  a $\C_p$-Banach space: it would be of dimension~$1$, 
 but we have a filtration
$$
\xymatrix@R=.4cm@C=.6cm{
0 \ar[r] & t\B^+_{\dr}/t^2\ar[r] \ar[d]^{\wr} & \B^+_{\dr}/t^2\ar[r] & \B^+_{\dr}/t\ar[r] \ar[d]^{\wr} &0\\
 & \C_p & & \C_p
 }
$$
whose graded pieces  are isomorphic to $\C_p$.
We can say that  $\B^+_{\dr}/t^2$ looks like a $\C_p$-vector space of dimension $2$; we will write this as $\B^+_{\dr}/t^2\sim \C_p^2$. 
More generally, we have $$\B^+_{\dr}/t^n\B^+_{\dr}\sim \C_p^n.$$
\end{example}
\subsubsection{de Rham and Hodge-Tate  comparison theorems}\label{comp-dr-theorem}
Define the period ring $\B_{\dr}:=\B^+_{\dr}[1/t]$. 
\begin{theorem}{\rm({\em de Rham comparison theorem}, Faltings, 1989 \cite{Fi})} Let 
$X$ be a proper and smooth variety over $K$, for $[K:\Q_p] <\infty$. There exists a period isomorphism
\begin{equation}
\label{comp-dr}
\alpha_{\dr}:\quad  H^n_{\dr}(X)\otimes_{K}\B_{\dr}\simeq H^n_{\eet}(X_{\ovk},\Q_p)\otimes_{\Q_p}\B_{\dr}
\end{equation}
compatible with the Galois action and filtration, where the Hodge filtration on de Rham cohomology is defined by
$$
F^iH^n_{\dr}(X) := {\rm Im}( H^n(X,\Omega^{\geq i}_{X})\to H^n_{\dr}(X) ),\quad i\geq 0.
$$
\end{theorem}
\begin{center}\fbox{
 $\B_{\dr}$ {\it  contains all $p$-adic periods of algebraic varieties over $K$}. 
 }
 \end{center}
\vskip.2cm

  The above theorem yields a filtered isomorphism (take $G_K$ fixed points of both sides of (\ref{comp-dr}) and use the fact that $\B_{\dr}^{G_K}=K$):
$$
H^n_{\dr}(X)\simeq (H^n_{\eet}(X_{\ovk},\Q_p)\otimes_{\Q_p}\B_{\dr})^{G_{K}}.
$$
However we are  not able to recover the Galois action on $H^n_{\eet}(X_{\ovk},\Q_p)$ from the left hand side of the period isomorphism (\ref{comp-dr}) (the structure both on de Rham cohomology and on $\B_{\dr}$ is too 
coarse 
for that). However, 
taking the $\gr^0_F$ of both sides of (\ref{comp-dr}), we obtain the following corollary 
(the conjecture of Tate~\cite{Ta} which was the starting point of all $p$-adic comparison theorems).
\begin{corollary}{\rm ({\em Hodge-Tate decomposition})}
 We have a Galois equivariant Hodge-Tate decomposition
\begin{equation}
\label{HT}
H^n_{\eet}(X_{\ovk},\Q_p)\otimes_{\Q_p}\C_p\simeq \bigoplus_{i\geq 0}H^{n-i}(X,\Omega^i_{X/K})\otimes_{K}\C_p(-i)
\end{equation}
\end{corollary}
This corollary say that twist of  the Galois representation $H^n_{\eet}(X_{\ovk},\Q_p)$ by $\C$ (with its natural Galois action) splits  as a direct sum of cyclotomic characters with multiplicities given by the Hodge numbers of the variety.
It is reminiscent of the classical Hodge decomposition for projective varieties over $\C$.
\subsection{Refined period rings} We will now introduce more sophisticated $p$-adic period rings, present the associated comparison theorems and   discuss their history. 
\subsubsection{The crystalline and semistable period rings} To recover the Galois representations on \'etale cohomology Fontaine defined more refined period rings:
$$
 \B_{\crr}\subset \B_{\st}\subset \ovk\cdot\B_{\st}\subset \B_{\dr}.
$$
They have the following properties:
\begin{enumerate}
\item The crystalline period ring 
   $\B_{\crr}$ is equipped with 
(commuting) 
$G_K$-action and Frobenius operator $\phi$. 
   \item The semistable period ring, besides $G_K$-action and Frobenius $\phi$,  has also an action of a $\B_{\crr}$-linear monodromy operator $N$ which is a $\B_{\crr}$-derivation (we have 
 $\B_{\st}^{N=0}=\B_{\crr}$) that commutes with the action of  $G_K$ and satisfies $N\varphi=p\varphi N$. 
   \end{enumerate}
\subsubsection{Comparison theorems} The most general $p$-adic comparison theorem for algebraic varieties takes the following form: 
   \begin{theorem}{\rm ({\em Semistable comparison theorem})} Let $X$ be a variety over $K$. 
   There exists a period isomorphism
   \begin{equation}
   \label{comp-st}
   \alpha_{\st}:\quad  H^n_{\hk}(X_{\ovk})\otimes_{K^{\nr}}\B_{\st}\simeq H^n_{\eet}(X_{\ovk},\Q_p)\otimes_{\Q_p}\B_{\st}, \quad n\geq 0,
   \end{equation}
   compatible with Frobenius, monodromy, Galois action,  and with the de Rham period isomorphism $\alpha_{\dr}$, i.e., $\alpha_{\st}\otimes\B_{\dr}\simeq \alpha_{\dr}$.
   \end{theorem}
   Here,   $H^n_{\hk}(X_{\ovk})$ is the Hyodo-Kato cohomology: a finite rank $K^{\nr}$-vector space\footnote{$K^{\nr}$ is the maximal unramified extension of the Witt vectors of the residue field  of $K$.}, equipped with a Frobenius, monodromy, and a locally constant
 action of $G_K$. Moreover, we have a Hyodo-Kato isomorphism
   $$
  \iota_{\hk}:  H^n_{\hk}(X_{\ovk})\otimes_{K^{\nr}}\ovk\simeq H^n_{\dr}(X_{\ovk}).
   $$
Hyodo-Kato cohomology plays the role of limit Hodge structures in the $p$-adic world. 
\vskip  3mm
 \begin{center}\fbox{
 $\ovk\cdot \B_{\st}$ {\it  contains  all $p$-adic periods of  algebraic varieties over $K$}. 
 }
 \end{center}
\vskip.2cm

 The period isomorphism (\ref{comp-st})  allows to go back and forth between de Rham cohomology and \'etale cohomology (with all the additional structures)
   \begin{align}\label{twist1}
   & H^n_{\eet}(X_{\ovk},\Q_p)\simeq (H^n_{\hk}(X_{\ovk})\otimes_{K^{\nr}}\B_{\st})^{N=0,\phi=1}\cap F^0(H^n_{\dr}(X_{\ovk})\otimes_{\ovk}\B_{\dr}),\\
   & H^n_{\hk}(X_{\ovk})\simeq (H^n_{\eet}(X_{\ovk},\Q_p)\otimes_{\Q_p}\B_{\st})^{G_K-{\rm sm}}.\notag
   \end{align}
   Here $(-)^{G_K-{\rm sm}}$ refers to taking vectors stable under the action of an open subgroup of $G_K$.

     \subsubsection{History of  $p$-adic Hodge Theory for algebraic varieties: comparison theorems}$\quad$

   (A) {\em  Early years: 1985-2010.}
     \begin{enumerate}
  \item  The first quite general proof of comparison theorems was given by Fontaine-Messing \cite{FM}. The proof worked for all large enough primes $p$ and also for some cohomologically simple varieties. It employed the syntomic technique, which was later refined by Hyodo and  Kato \cite{HK}, \cite{K1} (based on an earlier work of Bloch-Kato \cite{BK}), and finally  by Tsuji \cite{Ts} who proved the semistable comparison theorem in all generality. The modern version of the syntomic technique yields a computation of $p$-adic nearby cycles (which are, in general,   nontrivial for good reduction varieties) via syntomic cohomology defined as a filtered Frobenius eigenspace of crystalline cohomology. This coupled with the Poincar\'e duality or the rigidity supplied by Banach-Colmez spaces (see below) yields the comparison theorems for algebraic varieties. But, perhaps surprisingly, it has turned out also to be a powerful tool in studying comparison theorems for rigid analytic varieties (see below). So the technique is seeing a major revival recently. 
  \item  The first completely general proof of the de Rham\footnote{This was preceded by a proof of the Hodge-Tate decomposition \cite{Fa0}.} and the crystalline comparison theorems was given by Faltings \cite{Fi} who also  proved the semistable comparison theorem \cite{Fa} (a bit later than Tsuji). He uses the technique of almost \'etale extensions, which he invented for this purpose by refining the original ideas of Tate in dimension $0$. He is also able to treat local systems. This technique, amplified by a  strong almost purity result of Scholze~\cite{Sch0}, 
is very powerful and was also used recently to study comparison theorems for rigid analytic varieties (see below).
  \item  The author, in her study of Galois actions on geometric motivic cohomology, found yet another proof of the comparison theorems \cite{N4}, \cite{N2}: in it, both \'etale cohomology and syntomic cohomology are seen as   incarnations of $p$-adic motivic cohomology and the comparison isomorphism becomes an incarnation  of  the  localization map in motivic cohomology. The relation between motivic cohomology and $p$-adic cohomologies was studied extensively
% in the past 
by Geisser, Hesselholt, and Levine \cite{Hes}, \cite{GL}; a more recent study was initiated by Nikolaus-Scholze \cite{NSc} and Bhatt-Morrow-Scholze \cite{BMS2}.
  \end{enumerate}
  
  (B) {\em Recent years:  2010-}
  \begin{enumerate}
  \item   More recently,  Beilinson has found a new approach to comparison theorems \cite{BE1}, \cite{BE2} (Bhatt gave a similar proof of the semistable conjecture  in \cite{BH}). He uses $h$-topology (built from proper and Zariski open maps) to prove a Poincar\'e Lemma: the complex of differentials modulo $p^n$ becomes just constants in the $h$-topology; it is easy to see that mod $p^n$ differentials can be killed by $p$-covers; 
it is more difficult to kill  regular nonconstant functions mod $p^n$: for that Beilinson used a result of Bhatt's thesis.   Beilinson's approach is particularly simple in the case of the de Rham comparison theorem. 
  \item Scholze has rewritten the foundations of the almost \'etale approach to relative $p$-adic Hodge Theory \cite{Sch0}. He proved a general almost purity result (weaker versions were proved earlier by Faltings and Gabber-Ramero \cite{GR} and were a backbone of the almost \'etale  approach to $p$-adic comparison theorems) that allows to do $p$-adic Hodge Theory on algebraic varieties over $K$, and to extend it to rigid analytic varieties over $K$ (see below). It also allows to treat coefficients.
  \item  {\em Integral $p$-adic Hodge Theory.} We focus in this survey on rational $p$-adic Hodge Theory. But there exists  also integral $p$-adic Hodge Theory which is very important in applications to arithmetic. From the beginning of $p$-adic Hodge Theory we knew that we can do integral computations for large primes $p$ (Fontaine-Lafaille Theory) and also that the denominators that appear in general can be universally bounded 
(in terms of $p$ and the dimension of the varieties). 
However, perhaps no expert expected the result proved recently by Bhatt-Morrow-Scholze \cite{BMS}, \cite{BMS2} (see Bhatt-Scholze \cite{BS} for a different treatment via prismatic cohomology and \v{C}esnavi\v{c}ius-Koshikawa \cite{CK} for generalizations): it is enough to "untwist" the element $\mu\in\B_{\crr}$ (a lift of $\zeta_p-1$, for a primitive $p$'th root of unity $\zeta_p$,  from  $\C_p$)
 to obtain optimal integral $p$-adic comparison theorems. 
    \end{enumerate}
    
 (C)  {\em Uniqueness of $p$-adic comparison morphisms.} The comparison morphisms are normalized to be compatible with Chern classes (crystalline, de Rham, and \'etale). Hence one would expect (no matter the technique used for their construction) them to be all equal. The author proved in \cite{Ni}, \cite{Nii} this to be indeed the case for the syntomic, the old almost \'etale, and the Beilinson approaches;  the case of the (rational) comparison morphism of Bhatt-Morrow-Scholze and \v{C}esnavi\v{c}ius-Koshikawa was treated in the PhD thesis of Sally Gilles \cite{Gil}.
  The proof in \cite{Gil} is direct using the fact that the comparison morphism constructed in \cite{Gil} (a global geometric version of the morphism constructed by Colmez-Nizio{\l}~\cite{CN1}) 
is very similar to that of \cite{BMS} and \cite{CK} and can be shown directly to be the same as the comparison morphism of Fontaine-Messing. The proofs of the author exploit the motivic approach (A3) in \cite{Ni} and the $h$-topology approach (B1) in \cite{Nii}. These approaches allow to reduce all comparison morphisms to motivic localization morphism  or to the tautological morphism  in the Fundamental Exact Sequence of $p$-adic Hodge Theory (see below (\ref{seq1}))\footnote{The Fundamental Exact Sequence can be seen as an isomorphism between a Tate twist and a syntomic complex.}, respectively.

   \subsection{Banach-Colmez spaces}  Banach-Colmez spaces are omnipresent in modern $p$-adic Hodge Theory. We will illustrate how they enter the picture, discuss their properties, and give some examples. 
   \subsubsection{The fundamental exact sequences}
    Twisting the isomorphism (\ref{twist1}) by $\Q_p(r), r\geq 0$, we obtain the short exact sequence
   \begin{equation}\label{kwaku-kwak}
   0\to H^r_{\eet}(X_{\ovk},\Q_p(r))\to (H^r_{\hk}(X_{\ovk})\otimes_{K^{\nr}}\B^+_{\st})^{\N=0,\phi=p^r}\to (H^r_{\dr}(X_{\ovk})\otimes_{\ovk}\B^+_{\dr})/F^r\to 0
   \end{equation}
   Here $\B^+_{\st}$ is the ${\,}^+$-version of $\B_{\st}$, i.e., we have $\B_{\st}\simeq\B^+_{\st}[1/t]$.  We can ask: 
\vskip.2cm
  \begin{center}
  {\em In which category the above sequence lives ?}
  \end{center}
\vskip.2cm
The first term of  the sequence is a finite rank $\Q_p$-vector space, the last term, by Example \ref{derham}, looks like $\C_p^n$. What about the middle term ? Let us look at an example. 
\begin{example}Let $U=\{x=(x_0,x_1,\ldots),\  x_n\in B(1,1^{-}),\  x^p_{n+1}=x_n\}$. Consider the following commutative diagram
   $$
   \xymatrix{
 U\ar[d]^{x\mapsto x_0}\ar[r]^{\wt{\log}} & \B^+_{\dr}\ar[d]^{\theta}\\
  B(1,1^{-})\ar[r]^-{\log} & \C_p,
   }
   $$
   where  $\wt{\log}$ is a lifting of the logarithm such that 
$\wt{\log}\big((e^{\frac{2\pi i}{p^n}})_{n\in\N}\big)=t$. It induces  an isomorphism  $\wt{\log}: U\stackrel{\sim}{\to} \B_{\crr}^{+,\phi=p}.$
   Hence the above diagram and the exact sequence (\ref{logparis}) yield  the short exact sequence $$
   0\to \Q_pt\to \B^{+,\phi=p}_{\crr}\to \C_p\to 0
   $$
   We have obtained that  $\B^{+,\phi=p}_{\crr}\sim {\mathbf C}_p\oplus \Q_p$.
     
      More generally, for $m\geq 0$, we have the {\em Fundamental Exact Sequence of $p$-adic Hodge Theory}:
   \begin{equation}
   \label{seq1}
   0\to \Q_pt^m\to \B^{+,\phi=p^m}_{\crr}\to \B^+_{\dr}/t^m\B^+_{\dr}\to 0
   \end{equation}
   Hence  $ \B^{+,\phi=p^m}_{\crr}\sim {\mathbf C}_p^m\oplus \Q_p$.
   \end{example}
         
         Coming back to the exact sequence (\ref{kwaku-kwak}) we see that the terms are $\Q_p$- or $\C_p$-vector spaces, or extensions of such. The most naive guess for the category which contains such objects and in which one can do homological algebra would be the category of (locally convex) topological vector spaces over $\Q_p$. But this category is too ``flabby": there exists an  isomorphism 
${\mathbf C}_p\oplus\Q_p\simeq {\mathbf C}_p$ (to define this isomorphism consider both sides as $\Q_p$-Banach spaces and recall that in the $p$-adic world those behave like Hilbert spaces: they have orthonormal bases)
 and this we definitely do not want.

   \subsubsection{Banach-Colmez spaces}  We have however: 
   \begin{theorem}{\rm (Colmez, Fontaine, \cite{CB}, \cite{Fo-Cp})}
   There exists an abelian category of Banach-Colmez vector spaces $\mathbb{W}$ which are finite dimensional ${\mathbf C}_p$-vector spaces $\pm$ finite dimensional $\Q_p$-vector spaces. Moreover
   \begin{enumerate}
   \item Such a space has a total dimension 
${\rm Dim}(\mathbb{W}):=(\dim_{{\mathbf C}_p}\mathbb{W},\dim_{\Q_p}\mathbb{W})\in (\N,\Z)$,
   \item ${\rm Dim}(\mathbb{W})$ is additive on short exact sequences.
   \end{enumerate}
      \end{theorem}
For example ${\rm Dim}(\C_p)=(1,0)$ and ${\rm Dim}(\C_p\oplus\Q_p)=(1,1)$.
%The total dimension is an invariant of Banach-Colmez spaces. 
It follows that in this category  $\C_p$ can not be isomorphic to $\C_p\oplus \Q_p$.
 %if it were, looking at total dimensions, we would get an
     %equality $(1,0)=(1,1)$.  

\smallskip
 We will now discuss the above theorem in more detail.  
    A Vector Space (VS for short) ${\mathbb W}$ is a functor $\Lambda\to {\mathbb W}(\Lambda)$ from nice $\Q_p$-algebras to $\Q_p$-vector spaces. Here, a Banach algebra $\Lambda$ ($|x+y|\leq \sup (|x|, |y|)$, $ |xy|\leq |x||y|)$ is nice\footnote{Colmez called them sympathetique. They are perfectoid in the terminology of Scholze. } if $|x|=\sup |s(x)|_p, s:\Lambda\to \C_p$, and the map $x\mapsto x^p$ is surjective.
   \begin{example}
   \begin{enumerate}
   \item If $V$ is a finite dimensional $\Q_p$-vector space, 
the associated VS is the functor
$$\Lambda\mapsto V(\Lambda)=V,\quad (\Lambda_1\to \Lambda_2)\mapsto (\id: V\to V).$$ We have $V(\C_p)=V.$
   \item Let $d\geq 1$. The VS ${\mathbb V}^d$ is the functor 
$$ \Lambda\mapsto {\mathbb V}^d(\Lambda)=\Lambda^d,\quad (\Lambda_1\to \Lambda_2)\mapsto (\Lambda^d_1\to \Lambda^d_2).
   $$ We have ${\mathbb V}^d(\C_p)=\C_p^d$. 
   \end{enumerate}
   \end{example}
   A Vector Space ${\mathbb W}$ is called finite dimensional (nowadays called simply a Banach-Colmez space)  if it admits a presentation
   $$
   \xymatrix@R=.3cm@C=.5cm{
   0 \ar[r] & V_2\ar[rd]\\
   0\ar[r] & V_1\ar[r] & {\mathbb W}^{\prime}\ar[rd]\ar[r] & {\mathbb V}^d\ar[r] & 0\\
   &&& {\mathbb W}\ar[r]
   &  0
   }
   $$
   We have 
   ${\rm Dim}({\mathbb W})=(d, \dim_{\Q_p}V_1-\dim_{\Q_p}V_2)$. This value is independent of the presentation chosen (a  fact difficult to prove). The proof of this independence can also be used to show that
every Banach-Colmez space  ${\mathbb W}$ is a Banach Space 
(i.e., ${\mathbb W}(\Lambda)$ is a Banach space for any $\Lambda$).    
  \begin{example}We list some examples of Banach-Colmez spaces: 
   \begin{enumerate}
   \item $\B^+_{\dr}/t^m$ is $\mathbb{B}_m(\C_p)$ with ${\rm Dim} (\mathbb{B}_m)=(m,0).$
   \item $\B^{+,\phi^a=p^b}_{\crr}$ is $\mathbb{U}_{a,b}(\C_p)$ with ${\rm Dim}( \mathbb{U}_{a,b})=(b,a).$
 The exact sequence (\ref{seq1}) can be lifted to the category of Banach-Colmez spaces, i.e., it 
   is the $\C_p$-points of the exact sequence of Banach-Colmez spaces
   $$
   0\to \Q_pt^m\to {\mathbb U}_{1,m}\to {\mathbb B}_m\to 0 
   $$
   \item ${\mathbf C}_p/\Q_p$ is the $\C_p$-points of a Banach-Colmez space with ${\rm Dim}=(1,-1).$
   \end{enumerate}
   \end{example}
\begin{remark}
The category of Banach-Colmez spaces is quite rigid: If ${\mathbb W}_1$ is a successive extension of ${\mathbb V}^1$'s, and if ${\mathbb W}_2$ is of $\C_p$-dimension $0$, then any morphism $ {\mathbb W}_1\to {\mathbb W}_2$
is the $0$-map. This fact is very useful: it allows to circumvent Poincar\'e 
duality
arguments in the proofs of $p$-adic comparison theorems \cite{CN1}. 
\end{remark}
  \begin{remark}The idea of defining a category akin to the category of Banach-Colmez spaces goes back to Fontaine. 
Fontaine  studied similar  kind of structures equipped with Galois action (he called them ``almost $\C_p$-representations"). In Fontaine's category we can also distinguish between $\C_p$ and $\C_p\oplus\Q_p$ (eguipped with the natural action of $G_{K}$). Namely, if you believe that you can figure out the presence of  $\C_p$ in your space (because $\C_p$  is big) then, if you have $G_{K}$-action, you can take a look at the associated Euler characteristic of both spaces. It kills $\C_p$ and allows you to  figure out the presence of $\Q_p$. For Banach-Colmez spaces, Colmez found an analytic replacement for this Galois-theoretical argument.
  
  The theory of Banach-Colmez spaces was refined  and extended by Pl\^{u}t \cite{Plu}, and, more recently, by Fargues-Fontaine \cite{FF}, Fargues \cite{Far}, and Le-Bras \cite{LeB}. Banach-Colmez spaces  are a special case of Scholze's  diamonds \cite{SW}. 
  \end{remark}
\subsection{Complements}We will finish the survey of $p$-adic Hodge Theory for algebraic varieties with a brief discussion of its arithmetic aspects and number-theoretical applications.
\subsubsection{Arithmetic $p$-adic Hodge Theory} What we have discussed above is  the geometric side of $p$-adic Hodge Theory. But $p$-adic Hodge Theory  has also its arithmetic side. In it Fontaine defined  subcategories of $p$-adic Galois representations (i.e., continuous representations of $G_K$ on finite rank vector spaces over $\Q_p$):
\begin{equation}
\label{categories}
\Rep(G_K)\varsupsetneq \Rep_{\htt}(G_K)\varsupsetneq  \Rep_{\dr}(G_K)= \Rep_{\pst}(G_K)\varsupsetneq \Rep_{\rm geom}(G_K).
\end{equation}
The categories $\Rep_{\htt}(G_K), \Rep_{\dr}(G_K),$ and $\Rep_{\pst}(G_K)$ of Hodge-Tate, de Rham, and potentially semistable representation, respectively, consist of  $p$-adic representations that satisfy an abstract form of the Hodge-Tate decomposition (\ref{HT}), the  de Rham comparison isomorphism (\ref{comp-dr}), or the potentially semistable comparison isomorphism (\ref{comp-st}), respectively. The category $\Rep_{\rm geom}(G_K)$ is the category of representations coming from geometry, i.e., from subquotients of  $p$-adic \'etale cohomology of algebraic varieties over $K$. 

     To illustrate the strictness of the first two inclusions in (\ref{categories}) we will look at some extension groups. Let $\Ext^1_{*}(\Q_p,\Q_p(i))$, $*= G_K, {\rm HT}, \dr$, denote  the extension groups of $\Q_p$ by $\Q_p(i)$ in the respective categories. We have (the second result was  conjectured by Fontaine~\cite{Annals}
and proved by Bloch-Kato~\cite{BK1})
\begin{align*}
\Ext^1_{*}(\Q_p,\Q_p)=&\ \begin{cases} K\oplus \Q_p & \mbox{ if } *=G_K,\\
\Q_p & \mbox{ if } *={\rm HT},
\end{cases}\\
 \Ext^1_{*}(\Q_p,\Q_p(i) )=&\ \begin{cases} K & \mbox{ if } *={\rm HT},\\
0 & \mbox{ if } *={\rm \dr},
\end{cases},\quad{\text{for $i<0$.}}
\end{align*}
(For $i>0$, all $\Ext^1_{*}(\Q_p,\Q_p(i) )$ coincide.)
The fact that the inclusion $ \Rep_{\dr}(G_K)\supset \Rep_{\pst}(G_K)$ is actually an equality is a $p$-adic version of the potentially unipotent monodromy theorem (always true, by Grothendieck,  for  $\ell$-adic Galois representations, i.e., on vector spaces over  $\Q_{\ell}, \ell\neq p$). Berger \cite{Ber} reduced the proof of this theorem to a conjecture of Crew \cite{Cre} on $p$-adic differential equations which was then  immediately proved by Andr\'e \cite{An}, Mebkhout \cite{Me}, and Kedlaya \cite{Ked}.  
The last inclusion in (\ref{categories}) follows from Theorem \ref{comp-dr-theorem}. This inclusion is strict: we did not put any restrictions on eigenvalues of Frobenius (we note however that these restrictions could be added, see \cite{Vol} for an example) but, as shown in \cite{MV}, even an addition of such restrictions  would not  have been sufficient to force the equality. 

  \begin{example}The characterization of geometric representations as potentially semistable is very powerful. Very early on in the the development of $p$-adic Hodge Theory it allowed 
  to prove a conjecture of Shafarevich saying
that there are no Abelian varieties over $\Q$
 with good reduction everywhere~\cite{Fab,Ab}.
The argument goes as follows: assume that such an Abelian variety exists. One looks at its first $p$-adic geometric \'etale cohomology group.  By assumption, it is a crystalline representation (a version of potentially semistable representation abstracting the cohomological features of varieties with good reduction) corresponding to the first de Rham cohomology group. The Hodge weights of the latter give then an upper bound on the ramification of the largest  field whose Galois group acts trivially on $H^1_{\eet}$ with $\Z/p$ coefficients. The lower bound is supplied by Odlyzko. The two bounds superimpose giving a contradiction.
  \end{example}
  
   \begin{remark}
  The development of arithmetic
 $p$-adic Hodge Theory was followed by a quest for a $p$-adic local Langlands correspondence. 
This is a program which studies the relationship between $p$-adic Galois representations and $p$-adic unitary representations of $p$-adic reductive groups. It was initiated by Breuil \cite{Breu} and a major progress in dimension $2$ was made by Colmez \cite{ColK}.
\end{remark}
\subsubsection{Number-theoretical applications} \label{applications}
The inclusion $\Rep_{\rm geom}(G_K)\subset \Rep_{\pst}(G_K)$ 
makes sense 
in the global setting, that is over number fields, as well. 
 There, amazingly, it tends to be an equality. That is, 
\vskip.2cm
 \begin{center}\fbox{{\em in a global situation, de Rham implies  geometric.}}
 \end{center}
 \vskip 2mm
  This was stated as a conjecture: 
\begin{conjecture}{\rm (Fontaine-Mazur, \cite{FMaz}, \cite{FonA})}
Suppose that $\rho: G_{\Q}\to {\rm GL}(V)$ is an irreducible $p$-adic representation which is unramified at all but finitely many primes and  $\rho|G_{\Q_p}$ is de Rham. Then there is a smooth projective variety $X/\Q$ and integers $i,j$ such that $V$ is a subquotient of $H^i_{\eet}(X_{\overline \Q},\overline{\Q}_p(j))$.
\end{conjecture}
Combined with the conjectural global Langlands correspondence,
this conjecture implies that $2$-dimensional representations of $G_\Q$ satisfying the above conditions
come from modular forms.  This is exactly what Wiles~\cite{Wi} proved in 1994, in a special case, on his way to  proving Fermat Last Theorem (FLT); his proof makes heavy use of integral $p$-adic Hodge theory.

 The Fontaine-Mazur conjecture is basically known in dimension 2 by the work of  Emerton \cite{Em}, Kisin \cite{Kis}, and Lue Pan \cite{Pan}. 
\begin{example} ({\em Diophantine applications.}) Recall that Fermat Last Theorem  states that there are no nontrivial positive integral solutions to the equation
\begin{equation}
\label{Fermat} 
x^n+ y^n= z^n,\quad n\geq 3.
\end{equation}
This is an example of a Diophantine problem and $p$-adic Hodge Theory has been  used with great success to solve such problems. For example, before Wiles' prove of FLT we knew that there is  a finite number of positive integral solutions to the equation (\ref{Fermat}) because this equation can be seen as defining a plane curve of genus $>1$ and  Faltings proved the Mordell Conjecture in 1983  \cite{FaM} (see the wonderful sketch  of the proof by Bloch \cite{BlM}):  
\begin{conjecture}{\rm (Mordell,  1922)} A curve over a number field $K$ of genus $>1$ has only finitely many $K$-rational points.
\end{conjecture}
Faltings' proof uses as one of the tools a very early version of $p$-adic Hodge Theory. Let us sketch his argument.
In the first step, Faltings reduced proving finiteness of rational points  to proving  finiteness of  different mathematical objects:  isomorphism classes of principally polarized abelian varieties of dimension $d$  defined over $K$ with good reduction outside a given finite set of primes $S$. This problem, in turn,  he subdivided into two problems:  (a) finiteness of isogeny classes of Abelian varieties of dimension $d$ with good reduction outside $S$ (b) finiteness of isomorphism classes in every isogeny class. $p$-adic Hodge Theory is used in part (b) via Raynaud's
 theorem about Galois actions on points of group schemes \cite{Ray}. 

 For a modern, more powerful,  variation on Faltings' proof see the recent work of Lawrence-Venkatesh \cite{LV}. For a non-abelian version (studying iterated integrals instead of just integrals) see the work of Kim \cite{Kim} for theoretical aspects and \cite{Bal} for an application to  counting rational points on a particularly stubborn curve. 
%isogeny classes can be classified by the isomorphism classes of rational Tate modules of Abelian varieties (dually, their $H^1_{\eet}$) as Galois representations. These representations are semisimple hence it suffices to show that they give rise to only finite number of different trace functions on the Galois group. 
\end{example}
\begin{example}({\em Selmer groups.})  The extension groups $\Ext^1_{ \dr}(\Q_p,V)$ appearing above
are usually denoted by $H^1_g(G_K,V)$ (where $g$ stands for ``geometric'', a notation introduced by Bloch and Kato~\cite{BK1}).  
The fact that they are often strictly smaller than $H^1(G_K,V)$ is crucial for many questions and made them
omnipresent in modern Algebraic Number Theory under the name of 
{\em Selmer groups}. They have been generalized in two ways:
  
  (1) {\em geometrically, for varieties over local fields}: 
they can be defined for any algebraic variety over $K$ \cite{NN}, \cite{DN} under the name of 
{\em syntomic cohomology groups} and are used as an approximation of  
$p$-adic motivic cohomology (a refinement of $p$-adic \'etale cohomology capturing classes coming from geometry). 
  
  (2) {\em globally (over number fields)}: 
they can be globalized and extended to all $H^i$  (and not only $H^1$); 
see \cite{NDu} for a direct construction and \cite{Ven} 
for reinterpretations via derived Galois deformation rings. 
  
  Both generalizations come into play
 in the study of special values of (complex and $p$-adic) $L$-functions.  They enter in the computation
of the $p$-adic valuation of complex $L$-values, i.e.,~Tamagawa Numbers~\cite{BK1} and in the definition of
$p$-adic regulators that appear in  the $p$-adic $L$-values~\cite{BDR}. 
\end{example}
\section{Rigid analytic varieties} 
We now pass to the relatively new study of $p$-adic Hodge Theory for $p$-adic analytic varieties (called rigid analytic varieties since the seminal work of Tate~\cite{Tat}).
   \subsection{Rigid analytic varieties} 
Rigid analytic varieties have, in general, very large cohomology groups. Moreover, even locally, while $\ell$-adic \'etale cohomology, for $\ell\neq p$, tends to behave like de Rham cohomology, $p$-adic \'etale cohomology exhibits a very different behaviour. 
\subsubsection{Balls and annuli}
We will illustrate what we said above with a couple of examples. 
   \begin{example}
   \begin{enumerate}
   \item Let 
   ${\mathbb D}$ be the open unit disk in $\C_p$. We have
   \begin{align*}
 & H^1_{\dr}({\mathbb D})=0,\quad H^1_{\eet}({\mathbb D},\Q_{\ell})=0,\quad \ell\neq p;\\
& H^1_{\proeet}({\mathbb D},\Q_p)=\so({\mathbb D})/\C_p.
\end{align*}
        Pro-\'etale cohomology used here  is a version  of \'etale cohomology defined by Scholze \cite{Sch}, in which infinite (but inverse limits of finite)
 \'etale covers are allowed.  The cohomology group $H^1_{\proeet}({\mathbb D},\Q_p)$ is so big because of the 
existence of Artin-Schreier coverings: $y=x^p-x$ defines an \'etale covering of ${\mathbb A}^1_{{\mathbf F}_p}$ (we note that we have ${\rm d}y/{\rm d}x=-1$) and, hence, also of ${\mathbb D}$.
  \item Let $X=\{z\in\C_p,\ r<|z|<s\}$ be an open annulus in $\C_p$. We have 
 \begin{align*}
& H^1_{\dr}({X})\simeq\C_p\simeq\C_p<{\rm d}z/z>,\quad H^1_{\eet}({X},\Q_{\ell})\simeq \Q_{\ell},\quad  \ell\neq p;\\
 \mbox{ and an exact sequence }\\
& 
   0\to \so({X})/\C_p\to H^1_{\proeet}({X},\Q_p)\to \Q_p<{\rm d}z/z>\to 0
\end{align*}
\end{enumerate}
   \end{example}
\subsubsection{History of  $p$-adic Hodge Theory for rigid analytic  varieties.}
 \begin{enumerate}
 \item 1967: Tate~\cite{Ta} asked whether there is a Hodge Theory for $p$-adic analytic varieties.
 \item 2010: Scholze responded positively to Tate's question by developing $p$-adic Hodge Theory for  adic spaces \cite{Sch0}, \cite{Sch}. 
 \end{enumerate}

The subject is now 
 a very active area of research. 
For example, we have seen substantial progress in $p$-adic Hodge Theory 
for proper or Stein analytic varieties (the two endpoints of the spectrum).

  \subsection{Proper  rigid analytic varieties} The $p$-adic Hodge Theory for proper rigid analytic is very similar  to the $p$-adic Hodge Theory for algebraic varieties.  In particular, we have
  \begin{enumerate}
  \item $p$-adic geometric \'etale cohomology $H^n_{\eet}(-,\Q_p)$ is of  finite rank \cite{Sch}.
  \item there is a de Rham comparison theorem (with coefficients) obtained by Scholze via the almost \'etale technique amplified by  his  almost purity result \cite{Sch}. See \cite{Diao}, \cite{Poin} for generalizations. 
  \item there is a potentially semistable comparison theorem for smooth and ("almost") proper varieties obtained by Colmez-Nizio{\l} via the syntomic technique that works surprisingly well in the analytic context \cite{CN1}, \cite{CN4}.  
  \item there are integral comparison theorems of Bhatt-Morrow-Scholze, \v{C}esnavi\v{c}ius-Koshikawa, and Bhatt-Scholze \cite{BMS}, \cite{BMS2}, \cite{CK}, \cite{BS} that introduce a new cohomology ($\A_{\rm inf}$-chomology or prismatic cohomology) as a go-between for  \'etale and de Rham cohomologies.
    \end{enumerate}
  \subsection{Stein rigid analytic varieties} 
At the other end of the spectrum, one finds Stein varieties (analytic analogs of affine varieties).
The key property of Stein rigid varieties is the fact that coherent sheaves have no higher cohomology. 

\subsubsection{Some examples}
We will start with some examples of computations in the case the Stein space  has a semistable formal model over the ring of integers $\so_K$ 
  (in this case all the irreducible components of the special fiber of the formal model 
   are proper). 
   
\begin{example}
 
    (i)  {\em Rigid analytic affine space $\mathbb{A}^d_K$} (see \cite{CN2}):
   $$ 
   H^r_{\proeet}(\mathbb{A}^d_{{\mathbf C}_p},\Q_p(r))\simeq \Omega^{r-1}(\mathbb{A}^d_{{\mathbf C}_p})/\ker d,\quad r\geq 1,$$
   where $d$ is the de Rham differential. The pro-\'etale cohomology is not of finite rank. Its Banach-Colmez total dimension  ${\rm Dim}=(\infty,0)$.

    (ii) {\em  Torus $\mathbb{G}_{m,K}^d$ of dimension $d$ (and $r\geq 1$)}:
we have a short exact sequence of Banach-Colmez spaces
   $$
 0\to   \Omega^{r-1}(\mathbb{G}_{m,{\mathbf C}_p}^d)/\ker d\to H^r_{\proeet}(\mathbb{G}_{m,{\mathbf C}_p}^d,\Q_p(r))\to\bigwedge^r\Q_p^d\to 0.
   $$
Here the wedge product 
   $\bigwedge^r\Q_p^d=\oplus_{i_1<\cdots<i_r} (\dlog z_{i_1}\wedge\cdots\wedge \dlog z_{i_r})\Q_p$. The  total dimension of $H^r_{\proeet}(\mathbb{G}_{m,{\mathbf C}_p}^d,\Q_p(r))$ is
   ${\rm Dim}=(\infty, \binom{d}{r}).$
 
    (iii) {\em Drinfeld half-plane $\Omega_K:=\mathbb{P}_K\setminus\mathbb{P}(K)$}: we have a short exact sequence of Banach-Colmez spaces
   $$
    0\to   \so(\Omega_{{\mathbf C}_p})/\ker d\to H^1_{\proeet}(\Omega_{{\mathbf C}_p},\Q_p(1))\to {\rm Sp}(\Q_p)^*\to 0,
   $$where 
   ${\rm Sp}(\Q_p)=\scc^{\infty}(\mathbb{P}(K),\Q_p)/\Q_p$  is the  (smooth) Steinberg representation of ${\rm GL}_2(K)$. The total dimension of $H^1_{\proeet}(\Omega_{{\mathbf C}_p},\Q_p(1))$  is ${\rm Dim}=(\infty,\infty)$.
\end{example}
 
   \begin{remark}
   (i) We have a similar result~\cite{CDN3}   for $\Omega^d_K$, the  Drinfeld symmetric space of any dimension $d>1$.
   
   (ii)  We have similar results~\cite{CDN1} for \'etale coverings of $\Omega_K$.  
This is used to show that, for $K=\Q_p$, the $p$-adic \'etale cohomology of these coverings
encodes a part of the $p$-adic local Langlands correspondence, which yields a geometric realization for this correspondence (it is a classical result that the $\ell$-adic \'etale cohomology of these coverings can be used to provide a geometric realization of the classical local Langlands correspondence). 
   \end{remark}
\subsubsection{A comparion theorem}
  We will finish with the following general result:
   \begin{theorem}{\rm (Colmez-Dospinescu-Nizio{\l}, \cite{CDN3}, \cite{CN3}, \cite{CN4})} Let $r\geq 0$ and let  $X$ be a  Stein rigid analytic variety over $K$. There exists a $G_K$-equivariant      exact sequence:
   \begin{align*}
   0\to H^r_{\proeet}(X_{\C_p},\Q_p(r)) & \to \Omega^r(X_{\C_p})^{d=0}\oplus (H^r_{\hk}(X_{\C_p})\wh{\otimes}_{K^{\nr}}\B^+_{\st})^{N=0,\phi=p^r}\\
    & \verylomapr{\can+\iota_{\hk}} H^r_{\dr}(X_{\C_p})\to 0.
   \end{align*}
     \end{theorem}
     
      Hence, just as for algebraic varieties, in the world of Stein or  proper rigid analytic varieties 
\vskip.2cm
\begin{center}\fbox{ pro-\'etale cohomology can be recovered from de Rham data !} 
\end{center}

 \end{document}